\documentclass[jp]{imsart}

\RequirePackage[OT1]{fontenc}
\RequirePackage{amsthm,amsmath}
\RequirePackage[numbers]{natbib}
\RequirePackage[colorlinks,citecolor=blue,urlcolor=blue]{hyperref}

\startlocaldefs
\numberwithin{equation}{section}
\theoremstyle{plain}
\newtheorem{thm}{Theorem}[section]
\endlocaldefs

\begin{document}

\begin{frontmatter}
\title{Bounding the Maximum of Dependent 
Random Variables}
\runtitle{Bounding the Maximum}

\begin{aug}
\author{\fnms{J.A.} \snm{Hartigan}\ead[label=e1]{john.hartigan@yale.edu}}
\address{ Yale University\\
\printead{e1}}

\affiliation{Yale University}

\end{aug}

 \begin{abstract}
{\bf Abstract}: Let $M_n$ be the maximum of $n$ zero-mean  gaussian variables $X_1,..,X_n$ with covariance matrix of  minimum eigenvalue $\lambda$ and maximum eigenvalue $\Lambda$. Then, for $n \ge 70$,
$$\Pr\{M_n \ge \lambda \left (2 \log n - 2.5 - \log(2 \log n - 2.5) \right )^\frac{1}{2} -.68\Lambda\}  \ge \frac{1}{2}.$$
Bounds are also given for tail probabilities other than $\frac{1}{2}$. Upper bounds are  given for tail probabilities of the maximum of dependent identically distributed variables. As an application, the maximum of purely non-deterministic stationary Gaussian processes is shown to have the same first order asymptotic behaviour as the maximum of  independent gaussian processes.
\end{abstract}



\end{frontmatter}

\section{ Introduction. }

The asymptotic behaviour of the maximum $M_n$ of $n$ i.i.d.(independent and identically distributed) random variables with continuous distribution function $F$ is well known following Gumbel[2], namely that
 \[-\log n - \log \left (1-F(M_n) \right) \rightarrow G \mbox{ in distribution,}\] 
where the Gumbel variable $G = -\log( -\log U)$  for $U$ uniform. 

When the variables are not independent, but identically distributed, upper bounds for the tail probabilities of $M_n$ are available through 
\[ \Pr\{M_n \ge A\} \le  n\left(1-F(A)\right),\]
but lower bounds are rare. 

In [1], Berman examines conditions under which the maximum of a stationary process has a limiting distribution.

 If the process is exchangeable, then the variables are i.i.d conditional on some tail process, and the limiting maximum is a distributed as  mixture of  limiting distributions of iid maxima conditioned on the the tail process.  For example, if the process is gaussian, the limiting distribution will be that of $Y + M_n$ where $Y$ is gaussian, and $M_n$ is the maximum of $n$ i.i.d guassians independent of $Y$.

For stationary Gaussian processes with correlations satisfying
$$\lim_{|i-j| \rightarrow \infty}\rho(X_i, X_j)\log |i-j| =0,$$ the maximum behaves asymptotically like the maximum of i.i.d variables. Berman also presents conditions under which the maximum of a Markov chain has a limiting distribution.

 Resnick[3] shows that if the variables form a markov chain, the limiting distribution of the maximum may be reduced to the distribution of the maximum of a set of i.i.d random variables.

The principal result here  provides lower bounds for the tail probabilities of the maximum of dependent gaussian variables.\newline

\section{ General bounds}
We will use $\Pr X$ to denote the expectation of the random variable $X$, and $\{S\}$ to denote the function that is $1$ when $S$ is true, and $0$ when $S$ is false.
\begin{thm}
 Let $M_n$ denote  the maximum of n random variables $X_1, .. X_n$ each with continuous distribution function F . Then, for each $n$,  there exists an exponential variable $E$ with
\begin{equation}
 -\log n - \log \left (1-F(M_n) \right)  \le  E .
\end{equation}
\end{thm}
Proof:
Let $F_n$ denote the distribution function of $M_n$.  
 Then
\begin{eqnarray*}
 \{M_n > A\}  &\le& \sum_{i=1}^{n}\{X_i > A\},\\
 1-F_n(A) &\le&  n\left(1-F(A)\right),\\
-\log n -\log\left(1-F(M_n)\right) &\le& E,
\end{eqnarray*}
since $1-F_n(M_n) \sim \exp(-E)$.\newline

\begin{thm}
  Let $M_n$ denote  the maximum of n independent random variables each with continuous distribution function F . Then the function
$$ G(M_n) = -\log(- n\log F(M_n)), $$ has a Gumbel distribution and,
$$
G \le - \log \left (n(1-F(M_n)) \right) \le G + \exp(-G)/n.
$$
\end{thm}
Proof: Since $F(M_n)$ is the maximum of $n$ independent uniforms, $F(M_n) \sim U^\frac{1}{n}$, so
 $G = -\log(- n\log F(M_n))$ is Gumbel. Then

\begin{eqnarray*}
1 -F(M_n) &=& 1 - e^{-\frac{1}{n}\exp(-G)},\\
\frac{1}{n} \exp(-G)/(1 +  \frac{1}{n} \exp(-G)) &\le& 1 - e^{-\frac{1}{n}\exp(-G)} \le \frac{1}{n} \exp(-G),\\
G \le -\log \left (n(1-F(M_n))\right) &\le& G + \log(1 +  \frac{1}{n} \exp(-G)),\\
G \le -\log \left (n(1-F(M_n))\right) &\le& G + \frac{1}{n} \exp(-G).
\end{eqnarray*}
 
	It follows that the limiting distribution of $-\log n - \log \left (1-F(M_n) \right) $ is the Gumbel distribution. Note that $E$ and $G$ are very close in their tail distributions, so there is not much difference between the upper bounds in the independent and dependent cases. \newline

\section{ Gaussian bounds}

In the Gaussian case, we invert the standard tail bounds  for  $1-\Phi(x)$ for large $x$ so that we can accurately determine the asymptotic distribution of the maximum.
\begin{thm}
 Let $V=-2\log(1-\Phi(x))- \log(2\pi)$. 
\begin{eqnarray*}
\mbox{ For } x\ge 2, &V-\log V \le x^2 ,\\
\mbox{ For } x\ge 1, & x^2 \le  V - \log V + \log V/V .
\end{eqnarray*}
\end{thm}
Proof: The standard bounds from Abramowitz and Stegun(1972), p 932: \newline For $x \ge 1$, with $y=x^2$,
\begin{equation}
\begin{array}{rcl}
\frac{\phi(x)}{x}(1 -\frac{1}{x^2})\le & 1-\Phi(x)& \le\frac{\phi(x)}{x}(1-\frac{1}{x^2} + \frac{3}{x^4}),\\
y + \log y - 2\log\left(1 -\frac{1}{y} + \frac{3}{y^2}\right) \le& V& \le y + \log y  -2\log\left(1 -\frac{1}{y}\right).\\
\end{array}
\end{equation}

We demonstrate the specified bounds by explicit calculation for moderate $x$, and by using the standard bounds for large $x$. 
For the lower bound, let  $V-y-\log y= 2 \log \frac{dV}{dy}= \Delta$ where, from (3.1), 
$ \Delta \le  -2\log\left(1 -\frac{1}{y} )\right).$
 \begin{eqnarray*}
y-V+\log V &=& \log(V/y) - \Delta\\
&\ge& \log(1+\log y /y) +2\log\left(1 -\frac{1}{y} )\right) \\
&\ge& 0 \mbox{ for } y\ge 11.
\end{eqnarray*}
The lower bound is thus established for $y\ge 11$. The lower bound is exhibited in
 explicit calculation for $4<y<11$, so that the lower bound holds for $y \ge 4$ which is $x \ge 2$.

For the upper bound,  using $V=y+\log y + \Delta$ and $\Delta \ge \frac{1}{y}-\frac{3}{y^2}$,
 \begin{eqnarray*}
y-V+\log V -\log V/V &=& \log(V/y) - \log V/V - \Delta\\
&=& \log(1+(\log y +\Delta)/y) - \log V/V - \Delta\\
&\le& \log y/y - \log V/V  - (1-1/y)\Delta\\
&\le& (V-y)\log y/y^2 - (1-1/y)\Delta\\
&\le& (\log y +\Delta)\log y/y^2 - (1-1/y)\Delta\\
&\le& (\log y/y)^2 - (1-1/y-\log y/y^2)(\frac{1}{y}-\frac{3}{y^2})\\
 &\le& 0 \mbox{ for } y >5
\end{eqnarray*}

The upper bound is thus established for $y\ge 5$. The upper bound is exhibited in 
 explicit calculation for $1<y<5$, so that the upper bound holds for $x \ge 1$.\newline
\begin{thm}
 Let $M_n$ be the maximum of $n$ independent unit Gaussians. Let $N= \log(n^2/2\pi)$. For each $n$, there exists a Gumbel variable $G$, a monotone function of $M_n$, such that , for $M_n \ge 2$,

\begin{eqnarray}
(N+ 2G) - \log(2G+N) \le M_n^2 \le V - \log V + \log V/V
\end{eqnarray}
where $ V = N + 2G+ 2\exp(-G)/n$. 
\end{thm}

 Proof: Substitute the Gaussian probability bounds from lemma 1 into the probability bounds for $\Phi(M_n)$ given in Theorem 3.1. 

We see from theorem 3.2, that as $n \rightarrow \infty$,  $M_n^2 - N - \log N \rightarrow 2G$ in distribution. Note that the bounds apply only to tail probabilities $\Pr\{M_n \ge A\}$ for $A>2$.\newline

\begin{thm}
 Let $ M_n$ denote  the maximum of $n$ unit gaussian random variables. 
 Let $N= \log(n^2/2\pi)$. Then there is an exponential variable $E$ with
$$ M_n^2 \le \max \left(1,  N+2E - \log(N+2E) + \log(N+2E)/(N+2E)\right).$$
\end{thm}
Proof: From theorem 3.1, with $V= -\log(2\pi) - 2 \log(1-\Phi(M_n))$, we have, for $M_n \ge 1$,
$ M_n^2 \le V-\log V + \log V/V.$
 From theorem 1, there exists a function $E(M_n)$ distributed exponentially, with $V \le N+2E.$
Thus,
$$ M_n^2 \le \max \left(1,  N+2E - \log(N+2E) + \log(N+2E)/(N+2E)\right).$$

Asymptotically, $M_n^2 -N - \log N \le 2E$ as $n \rightarrow \infty$.\newline

\begin{thm}
 Let $M_n$ be the maximum of $n$ zero-mean  gaussian variables $X_1,..,X_n$.\newline
 Define
\begin{eqnarray*}
E_i &=& \Pr(X_i|X_1,..,X_{i-1}),\\
 R_i &=& X_i-E_i,\\
\tau^2 &=& \max_{1\le i\le n} \tau_i^2= \Pr E_i^2,\\ 
\sigma^2 &=& \min_{1\le i\le n} \sigma_i^2= \Pr R_i^2,\\ 
N&=& \log(n^2/2\pi),\\
L_\alpha &=& -2\log (-\log \alpha).
\end{eqnarray*}

Then, for $N + L_\alpha \ge 6$,   
$$\Pr\{M_n \ge \sigma(N+L_\alpha -\log (N+L_\alpha)^\frac{1}{2}+ \tau \Phi^{-1}(\alpha)  \} \ge 1-2\alpha.$$ 
 \end{thm}

Proof:  We first show that for each real $A$ and non-positive $B$,
$$ \Pr\{M_n \ge A + B\} \ge \Pr\{ \max_iR_i\ge A\} \min_i \Pr \{E_i \ge B\}.$$ \newline
Construct $n$ disjoint events
$$ H_i =\{R_i \ge A\}\prod_{j >i}\{R_j < A\}.$$
Note that $\sum_i H_i = \{ \max_i R_i \ge A \}$.
Since the terms $\{R_i \ge A\}$ are independent of all variables $X_j, j<i$,  the individual terms in each $H_i$ are l independent, and  $H_i$ is independent of $\{E_i \ge B\}$.
Also, $\{M_n \ge A + B\} \ge \sum_i H_i \{E_i \ge B\}$, since at most one of the events $H_i \{E_i \ge B\}$ occurs, and if any occurs $M_n \ge A+B$. Then

\begin{eqnarray*}
\Pr\{M_n \ge A + B\} &\ge& \sum_i \Pr H_i \Pr\{E_i \ge B\}\\
&\ge& \sum_i \Pr H_i \min_j\Pr\{E_j \ge B\}\\
&\ge&  \Pr \{ \max_i R_i \ge A \} \min_i\Pr\{E_i \ge B\}.
\end{eqnarray*}

For $A>0$, $\{\max_i R_i  \ge  A\} \ge \{\sigma \max_i (R_i /\sigma_i) \ge A\}$.

 Let $ R_{(n)} = \sigma \max_i (R_i /\sigma_i)$. From theorem 4, for $A/\sigma  \ge 2$, 
$$ \Pr \{ R_{(n)} \ge A\} \ge \Pr\{\left ( N+2G - \log ( N+2G)\right )^\frac{1}{2} \ge A/\sigma\}. $$
Since $\Pr \{2G \ge L_\alpha\} = 1-\alpha$, for $N + L_\alpha \ge 6$ (which implies  $A/\sigma  \ge 2$),
$$ \Pr \{ R_{(n)} \ge  \sigma (N+L_\alpha -\log (N+L_\alpha))^\frac{1}{2}\} \ge 1-\alpha.$$
Also $\min_i \Pr \{E_i \ge B\} = 1-\Phi(B/\tau$), so
$$\min_i \Pr \{E_i \ge \tau \Phi^{-1}(\alpha)\} = 1-\alpha.$$
Combining the two bounds gives, for $N + L_\alpha \ge 6$,
 $$\Pr\{M_n \ge \sigma(N+L_\alpha -\log (N+L_\alpha)^\frac{1}{2}+ \tau \Phi^{-1}(\alpha)  \} \ge 1-2\alpha.$$ 
 as asserted.\newline

If the random variables $X_1,..,X_n$ have covariance matrix $C$, then
 \[\sigma^2 = \min_i \frac{1}{C_{ii}^{-1}} \ge \min_x\frac{x'Cx}{x'x},\]
 $$\tau^2 = \max_i \tau_i^2 \le \max_i \Pr X_i^2 \le \max_x \frac{x'Cx}{x'x}.$$
Thus the stated lower bound holds if the minimum eigenvalue of $C$ is substituted for $\sigma^2$ and the maximum eigenvalue of $C$ is substituted for $\tau^2$.

In simulations, the lower bound is reasonably close to the actual distribution of $M_n$  in the tail when the variables are i.i.d, but can be quite conservative when the minimum eigenvalue is small.

\section{Application to stationary Gaussian processes}

Following Wold[4], purely non-deterministic stationary Gaussian processes \newline $X_i,  \infty < i < \infty, $ may be expressed in terms of the i.i.d innovations $$Z_i = X_i -\Pr(X_i|X_{i-1}, X_{i-2},..):$$
\[X_i =Z_i +  \sum_{j=1}^{\infty}\psi_jZ_{i-j}.\]  

Let $X_i$ have variance 1. Theorem 5 now applies with $\sigma^2 = \Pr Z_0^2, \tau^2 = 1- \sigma^2.$  The inequality is dominated asymptotically by the first term, so that, 
$$\mbox{ for each } \epsilon > 0, \Pr\{M_n > \sqrt{ 2\sigma^2 \log n }(1-\epsilon)\} \rightarrow 0 \mbox{ as } n \rightarrow \infty.$$

 The subsampled series $Y_i = X_{ik}$ for integer $k$ has innovations $Z_i +\sum_{j=1}^{k-1}\psi_jZ_{i-j}$ with residual variances $1 - \sigma^2(\sum_{j=k}^\infty \psi_j^2)$, which may be chosen arbitrarily close to $1$ by choosing $k$ large enough. Since $M_n = \max_1^{n} X_i \ge \max_1^{n/k} Y_i$, the stated inequality for $M_n$ holds with these residual variances. \newline
Thus, for each $\epsilon > 0$,  $M_n > \sqrt{2 \log n}(1-\epsilon)$ in probability as $n \rightarrow \infty$. (The order of magnitude based on $Y$ is   $\sqrt{2 \log (n/k)}$, but the $k$ washes out asymptotically.) From theorem 3.2, we have
$M_n < \sqrt{2 \log n}(1+\epsilon)$  in probability.  Thus $M_n/\sqrt{2 \log n} \rightarrow 1$ in probability, just as for an i.i.d. series. Detailed lower bounds for the tail probabilities of $M_n$ may be obtained depending on the rate of convergence of the series $\sum_{j=1}^\infty \psi_j^2$.


\begin{thebibliography}{9}
\bibitem{r1}
\textsc{Abramowitz, M.}and \textsc{Stegun, I. A., eds.} (1972).
\textit{Handbook of Mathematical Functions with Formulas, Graphs, and Mathematical Tables}.
Dover, New York MR0208797

\bibitem{r2}
\textsc{Berman, S. I.}(1964).
Limit Theorems for the Maximum Term in Stationary Sequences
\textit{Ann. Math. Statist.}
\textbf{35} 502-516. MR0161365

\bibitem{r3}
\textsc{Gumbel, E. J.}(1941)
The Return Period of Flood Flows
\textit{Ann. Math. Statist.}
\textbf{12} 163-190. MR0004457

\bibitem{r4}
\textsc{Resnick, S. I.}(1972)
Stability of Maxima of Random Variables Defined on a Markov Chain
\textit{Adv. Appl. Prob.}
\textbf{4} 284-295. MR0336820

\bibitem{r5}
\textsc{Wold, H. A.}(1938)
\textit{ Study in the Analysis of Stationary Time Series.}
 Almqvist and Wiksel, Uppsala

\end{thebibliography}
\end{document}